\documentclass[A4paper,reqno,oneside,12pt]{amsart}
\usepackage{fullpage}
\usepackage{setspace}
\usepackage{datetime}
\usepackage{tikz}
\usepackage{mathtools}
\usetikzlibrary{er}
\usetikzlibrary{arrows.meta}
\usepackage{amsmath,amsfonts,amssymb}
\usepackage{verbatim,enumitem,graphics}
\usepackage{microtype}
\usepackage{xcolor}
\usepackage[backref=page]{hyperref}
\hypersetup{%
    colorlinks,
    linkcolor={red!50!black},
    citecolor={green!50!black},
    urlcolor={blue!50!black}
}
\usepackage{dsfont}
\usepackage[abbrev,msc-links,backrefs]{amsrefs} 
\usepackage{doi}
\usepackage{lineno}
\usepackage[normalem]{ulem}

\newtheorem{theorem}{Theorem}[section]

\newtheorem{proposition}[theorem]{Proposition}

\newtheorem{problem}[theorem]{Problem}

\theoremstyle{definition}

\newcommand{\N}{\mathbb{N}}

\newcommand{\bbE}{\mathbb{E}}

\allowdisplaybreaks

\begin{document}

\setstretch{1.27}

\title{Sharp exponents for bipartite Erd\H{o}s-Rado numbers}

\author{D\'aniel Dob\'ak}

\author{Eion Mulrenin}


\address{Department of Mathematics, Emory University, Atlanta, GA, 30322, USA}

\email{\{daniel.dobak|eion.mulrenin\}@emory.edu}


\begin{abstract}
    The Erd\H{o}s-Rado canonization theorem generalizes Ramsey's theorem to edge-colorings with an unbounded number of colors, in the sense that for $n = ER(m)$ sufficiently large, any edge-coloring of $E(K_n) \to \N$ will yield some copy of $K_m$ which is colored according to one of four canonical patterns.
    In this 
    paper,
    we show that in the bipartite setting, the bipartite Erd\H{o}s-Rado number $ER_B(m)$ satisfies
    \begin{equation*}
        \log ER_B(m) = \Theta(m \log m).
    \end{equation*}
    Comparing this to the non-bipartite setting, the best known lower and upper bounds on $\log ER(m)$ are still separated by a factor of $\log m$.
\end{abstract}

\maketitle

\section{Introduction}

\subsection{Canonical Ramsey theory.}
The 
Erd\H{o}s-Rado canonization theorem \cite{ErdosRado2}, first proved by its eponymous authors in the 1950's, extends Ramsey's theorem \cite{Ramsey} to colorings with an arbitrary number of 
colors.
In contrast to Ramsey's theorem, no longer can one hope to ask for a monochromatic clique under any edge-coloring; for example, it may be the case that each edge is assigned its own unique color.
It is also impossible to demand just a monochromatic or rainbow clique, as is demonstrated by the so-called \textit{lexical} colorings in which one orders the vertex set and uses unique colors for all edges sharing the same minimum (or maximum) vertex.
The content of the 
canonization theorem is that, if the size of the host graph is large enough, some clique must be colored with one of these four ``canonical" color patterns.

\begin{theorem}[Erd\H{o}s-Rado canonization theorem for graphs]
    For all positive integers $m \geq 2$, there exists a least integer $n = ER(m)$ such that for any coloring\footnote{Here and throughout, we will use the following standard notation: for a positive integer $n$, $[n] = \{1, \dots, n\}$; and for a set $X$ and a non-negative integer $k$, $X^{(k)} = \{A \subseteq X: |A|=k\}$.} $\Delta: [n]^{(2)} \to \N$, there is a set $M \subseteq [n]$ of size $|M| = m$ with the following property: for all pairs $\{u_1, v_1\}, \{u_2, v_2\} \in M^{(2)}$ with $u_i < v_i$, either:
    \begin{itemize}
        \item $\Delta(\{u_1,v_1\}) = \Delta(\{u_2,v_2\})$ always holds;
        \item $\Delta(\{u_1,v_1\}) = \Delta(\{u_2,v_2\})$ holds if and only if $u_1 = u_2$;
        \item $\Delta(\{u_1,v_1\}) = \Delta(\{u_2,v_2\})$ holds if and only if $v_1 = v_2$;
        \item $\Delta(\{u_1,v_1\}) = \Delta(\{u_2,v_2\})$ holds if and only if both $u_1 = u_2$ and $v_1 = v_2$.
    \end{itemize}
\end{theorem}

\noindent The four coloring patterns of $M^{(2)}$ above are called, respectively, \textit{monochromatic}, \textit{min-colored}, \textit{max-colored}, and \textit{rainbow}.

In comparison with the quantitative bounds for the usual Ramsey numbers $R(m)$, which have been known to satisfy $\log R(m) = \Theta(m)$ since the 1940's, the situation for canonical Ramsey numbers is much worse.
One may observe that since the three canonical colorings which are not monochromatic each require at least $m-1$ colors, an easy lower bound for Erd\H{o}s-Rado numbers comes directly from the lower bound for the $(m-2)$-color Ramsey number $R(m; m-2)$ (see, e.g., \cites{Abbott, ConlonFerber}):
\begin{equation*}
    ER(m) \geq R(m; m-2) \geq 2^{c \cdot m^2},
\end{equation*}
where $c > 0$ is an absolute constant.
On the other hand, in the 1990's, Lefmann and R\"odl \cite{LefmannRodl} obtained the upper bound
\begin{equation*}
    ER(m) \leq 2^{c' \cdot m^2 \log m}
\end{equation*}
for an absolute constant $c' > 0$, improving on the triple-exponential bound given by the original proof of Erd\H{o}s and Rado. 
Therefore, the current state of the art for Erd\H{o}s-Rado numbers\footnote{We remark, however, that, unlike for 2-color Ramsey numbers, the tower height for hypergraph Erd\H{o}s-Rado numbers is known from the work of Lefmann and R\"odl \cite{LefmannRodl} and Shelah \cite{Shelah2} (and see \cite{Quantitative} for a further improvement to the upper bound).} is
\begin{equation}\label{canbounds}
    c \cdot m^2 \leq \log ER(m) \leq c' \cdot m^2 \log m.
\end{equation}
We will expand on where the logarithm in the upper bound comes from in Section \ref{conclusion}.
For a more detailed background discussion, see the recent paper \cite{KamcevSchacht}.

\subsection{The bipartite canonization theorem.}
Several years after Erd\H{o}s and Rado's seminal paper was first published, Rado \cite{Rado} obtained (as a corollary of much more general results | cf. \cite{Voigt} and \cite{Promel}, \S 9.3) a bipartite version of the canonization theorem.

\begin{theorem}[Bipartite canonization theorem]
    For all positive integers $m$, there is a least integer $n = ER_B(m)$ such that for any coloring $\Delta: E(K_{n,n}) \to \N$, there exists a copy of $K_{m,m}$ with the property that for all of its edges $(a_1, b_1), (a_2, b_2)$, either:
    \begin{itemize}
        \item $\Delta(a_1,b_1) = \Delta(a_2,b_2)$ always holds;
        \item $\Delta(a_1,b_1) = \Delta(a_2,b_2)$ holds if and only if $a_1 = a_2$;
        \item $\Delta(a_1,b_1) = \Delta(a_2,b_2)$ holds if and only if $b_1 = b_2$;
        \item $\Delta(a_1,b_1) = \Delta(a_2,b_2)$ holds if and only if both $a_1 = a_2$ and $b_1 = b_2$.
    \end{itemize}
\end{theorem}

\noindent We will refer to these four color patterns, and the colorings of $K_{m,m}$ thereof, respectively, as \textit{monochromatic}, \textit{left-colored}, \textit{right-colored}, and \textit{rainbow}.

The previously known quantitative bounds for $ER_B(m)$ were decently far apart.
The lower bound, coming from a random coloring, is exponential in $m \log m$, whereas the upper bound from Rado's original proof is triple-exponential in $m$.
Our main result determines precisely the order of magnitude of $\log ER_B(m)$ up to constant terms, in contrast to the gap persisting in the bounds for $\log ER(m)$ in (\ref{canbounds}).

\begin{theorem}\label{main}
    There exist absolute constants $c, c' > 0$ such that for all positive integers $m$,
    \begin{equation*}
        2^{c \cdot m \log m} \leq ER_B(m) \leq 2^{c' \cdot m \log m}.
    \end{equation*}
\end{theorem}

Our proof shows that for large $m$, one may take the logarithm to base 2 and $c = 1 - o(1)$, $c' = 8 + o(1)$, where the $o(1)$ terms goes to zero as $m$ grows.
It would certainly be interesting to narrow the gap a bit more.

It is worth pointing out that one may adapt the proof of Lefmann and R\"odl to the bipartite setting to obtain an upper bound which is exponential in $m^2 \log m$, but due to the neighborhood-chasing component of their argument, this seems to be a quantitative barrier for their method even in this setting.
Our proof diverts from the Lefmann-R\"odl approach from the start by looking at the colorings vertices in the second set induce on the first set and then splitting into two regimes depending on whether many vertices induce colorings with monochromatic $m$-sets or not.
By starting with a sufficiently asymmetric complete bipartite graph, a simple pigeonhole argument in the first regime gives either a monochromatic or right-colored $K_{m,m}$.
In the other regime, we use the asymmetric version of the K\H{o}v\'ari-S\'os-Tur\'an theorem to avoid a neighborhood-chasing approach.

Coincidentally, on the same day that the first version of this article was made public, a preprint by Gishboliner, Milojevi\'c, Sudakov, and Wigderson~\cite{GMSY24} was also posted in which the authors studied Erd\H{o}s-Rado numbers for sparse graphs (where the host graph is still a clique).
Among several other results, they show that every $n$-vertex bipartite graph $H$ with average degree $d$ has $ER(H) = n^{\Theta(d)}$.
Note that a special case of their result states that $ER(K_{m,m}) = 2^{\Theta(m \log m)}$, which may be compared with Theorem~\ref{main}.
We refer the interested reader to their paper for more details.

The remainder of our paper is organized as follows: in Section \ref{lowerbound}, we include a proof of the lower bound for 
completeness; in Section \ref{upperbound}, we give our proof of the upper bound; and in Section \ref{conclusion}, we offer some concluding remarks and discuss further directions and open problems.


\section{The lower bound in Theorem \ref{main}}\label{lowerbound}

\begin{proof}[\unskip\nopunct]
    Let $n=ER_B(m)$, and color the edges of $K_{n,n}$ with $m^2-1$ colors independently and uniformly at random. 
    Then there is no rainbow copy of $K_{m,m}$ since it would require at least $m^2$ colors, and so by our choice of $n$ it must contain a copy of $K_{m,m}$ that is either monochromatic, left-colored, or right-colored. 
    This implies that 
    \begin{align} \label{expeqn}
        \bbE[\text{number copies of }K_{m,m} \text{ that are either monochromatic or left- or right-colored}] \geq 1;
    \end{align}
    but also,
    \begin{align*}
        \bbE[\text{number of monochromatic copies of }K_{m,m}] &=\binom{n}{m}^2 \cdot \left(\frac{1}{m^2-1}\right)^{m^2-1} \\
        &\leq \left(\frac{en}{m}\right)^{2m} \cdot \left(\frac{1}{m^2-1}\right)^{m^2-1}
    \end{align*}
    and 
    \begin{align*}
        \bbE[\text{number of left-colored copies of }K_{m,m}] &\leq \binom{n}{m}^2 \cdot \left(\frac{1}{m^2-1}\right)^{m^2-m} \\
        &\leq \left(\frac{en}{m}\right)^{2m} \cdot \left(\frac{1}{m^2-1}\right)^{m^2-m}.
    \end{align*}
    Note that by symmetry the expected number of left- and right-colored copies
    of $K_{m,m}$ is equal. Using this together with (\ref{expeqn}), we obtain 
    \begin{align*}
        3\cdot \left(\frac{en}{m}\right)^{2m} \cdot \left(\frac{1}{m^2-1}\right)^{m^2-m} \ge 1
    \end{align*}
    from which 
    \begin{align*}
        n \ge \frac{m\cdot (m^2-1)^{\frac{m-1}{2}}}{3^{\frac{1}{2m}}\cdot e} = 2^{(1-o(1))m\log m}
    \end{align*}
    as claimed. 
\end{proof}


\section{The upper bound in Theorem \ref{main}}\label{upperbound}

We will use throughout the proof the following folklore bound for coloring singletons with arbitrarily many colors, sometimes referred to as the \textit{canonical pigeonhole principle}.

\begin{proposition}\label{CPHP}
    Let $ER_1(m)$ be the least integer $n$ such that any coloring $\Delta: [n] \to \N$ yields a set $M \subseteq [n]$ of size $m$ with $\Delta$ either constant, i.e., monochromatic, or one-to-one, i.e., rainbow, on $M$.
    Then
    \begin{equation*}
        ER_1(m) = (m-1)^2 + 1.
    \end{equation*}
\end{proposition}

\vspace{0.4cm}

\begin{proof}[\unskip\nopunct]
    Now, to begin, let $\Delta: E(K_{n_1,n_2}) \to \N$ be a given coloring, with 
    \begin{equation*}
        n_1 = 25 m^{9} \text{ and } n_2 = 2 \cdot \binom{n_1}{m} \cdot m^2.
    \end{equation*}
    By using $\binom{n_1}{m} \leq \left( \frac{e n_1}{m} \right)^m$, we have $n_1, n_2\leq 2^{(1+o(1)) 8 m \log m}$.

    Note that each $b\in [n_2]$ induces a coloring $\Delta_b$ of $[n_1]$ given by
    \begin{equation*}
        \Delta_b(a) = \Delta(a,b).
    \end{equation*}
    First, suppose there is a set $T_2 \subseteq [n_2]$ with $|T_2| \ge \binom{n_1}{m} \cdot m^2$ such that for each $b \in T_2$, $\Delta_b$ colors an $m$-set in $[n_1]$ monochromatically.
    Then by the pigeonhole principle, we can find $m^2$ such elements which all color the same $m$-set $T_1 \subseteq [n_1]$ monochromatically.
    Since $m^2 \geq ER_1(m)$, if we partition this $m^2$-set by the color each vertex uses on $T_1$, then we can find either a monochromatic or rainbow subset of size $m$, which will in turn yield a monochromatic or 
    right-colored
    $K_{m,m}$, respectively.

    Hence, assume that no such set $T_2$ exists and that instead there exists a subset $T_2' \subseteq [n_2]$ of size  $\binom{n_1}{m} \cdot m^2$ such that for all $b \in T_2'$, the coloring $\Delta_b$ contains no monochromatic $m$-set in $[n_1]$.
    In other words, for each $\Delta_b$, its color classes in $[n_1]$ all have size at most $m-1$. 
    Our goal is to find a subset $S_1 \subseteq [n_1]$ of size $4m^4+m$ and a subset $S_2 \subseteq T_2'$ of size $4^m \cdot m^{3m+1}$ (the reason for these choices will be clear later) such that for each $b \in S_2$ the induced coloring $\Delta_b$ is rainbow on $S_1$.

    Pick a $(4m^4+m)$-tuple $A = (a_1, \dots, a_{4m^4+m})$ from $[n_1]$ (possibly with repetitions) uniformly at random. 
    It suffices to show that
    \begin{equation*}
        \bbE\left[|\{ b \in T_2' : \Delta_b \text{ colors }A\text{ rainbow} \}| \right] \geq 4^m \cdot m^{3m+1},
    \end{equation*}
    as then we could just pick an $A$ for which this holds\footnote{Note that since the elements of $A$ receive different colors, we cannot have repetitions.} and take $S_1=A$ and $S_2$ some subset of $\{b \in T_2' :\Delta_b\text{ colors } A \text{ rainbow}\}$ of the desired size.
    
    To that end, fix a vertex $b \in T_2'$. 
    Since each color class of $\Delta_b$ contains at most $m-1$ vertices of $[n_1]$, the probability that any $a_i$ receives a certain fixed color is at most $\frac{m-1}{25m^9}$; and since the selection of the vertices $a_i$, $a_j$, $i < j$, is independent, we have
    \begin{align*}
        \mathbb{P}\left[\Delta_b(a_j) = \Delta_b(a_i) \right]
        &= \sum_{c \, \in \Delta_b([n_1])} \mathbb{P}\left[\Delta_b(a_j) = c \right] \cdot \mathbb{P} \left[ \Delta_b(a_i) = c \right] \\
        &\leq \frac{m-1}{25m^9} \cdot \sum_{c \, \in \Delta_b([n_1])} \mathbb{P} \left[ \Delta_b(a_i) = c \right]\\
        &= \frac{m-1}{25m^9}.
    \end{align*}
    Hence, by the union bound we have 
    \begin{align*}
        \mathbb{P}\left[\Delta_b\text{ does not color }A\text{ rainbow}\right]
        &= \mathbb{P}\left[ \bigcup_{1 \leq i < j \leq 4m^4+m} \Delta_b(a_i) = \Delta_b(a_j) \right]\\
        &\leq \sum_{1 \leq i < j \leq 4m^4+m} \mathbb{P}\left[\Delta_b(a_j) = \Delta_b(a_i) \right]\\
        &\leq \sum_{1 \leq i < j \leq 4m^4+m} \frac{m-1}{25m^9}\\
        &= \binom{4m^4 + m}{2} \cdot \frac{m-1}{25m^9}\\
        &\leq \frac{(5m^4)^2}{2} \cdot \frac{m-1}{25m^9} \\
        &\leq \frac{1}{2}.
    \end{align*}
    Therefore, $\mathbb{P}\left[\Delta_b\text{ colors }A\text{ rainbow}\right] \geq 1/2$, and so
    \begin{align*}
        \bbE\left[|\{ b \in T_2' : \Delta_b \text{ colors } A \text{ rainbow} \} | \right] &= \sum_{b \in T_2'} \mathbb{P}\left[\Delta_b\text{ colors }A\text{ rainbow}\right] \\ 
        &\geq \frac{1}{2} \cdot \binom{25m^9}{m} \cdot m^2 \\ 
        &\geq 4^m \cdot m^{3m+1}.
    \end{align*}

    Hence, we
    can choose $S_1 \subseteq [n_1]$ and $S_2 \subseteq [n_2]$ of sizes $|S_1| = 4m^4 + m$ and $|S_2| = 4^m \cdot m^{3m+1}$ such that for each $b \in S_2$, $\Delta_b$ colors $S_1$ with a rainbow pattern.
    It remains to find a canonically colored copy of $K_{m,m}$ inside the graph on $S_1 \times S_2$.
    For that, we will use the K\H{o}v\'ari-S\'os-Tur\'an theorem.
    
    \begin{theorem}[K\H{o}v\'ari-S\'os-Tur\'an, \cite{KST}]\label{KST}
        Let $G$ be a bipartite graph on $S_1 \times S_2$ which contains no copy of $K_{m,m}$.
        Then
        \begin{equation*}
            |E(G)| \leq (m-1)^{1/m} (|S_1|-m+1) |S_2|^{1-1/m} + (m-1) |S_2|.
        \end{equation*}
    \end{theorem}

    \textbf{Case 1.} Suppose that some set $M_1$ of $4m^4$ vertices in $S_1$ each have some ``popular" color, that is, a color which is used by at least $2 m^{1/m} \cdot |S_2|^{1-1/m}$ of their edges into $S_2$.
    Remove all edges from the graph which are not both incident to a vertex in $M_1$ and colored with its popular color, and remove all vertices from $S_1 \setminus M_1$.
    We are left with a bipartite graph $G$ on $M_1 \times S_2$ having $|M_1| = 4m^4$ and $|S_2| = 4^m \cdot m^{3m+1}$ and in which the number of edges is at least
    \begin{align*}
        8 m^{1/m} \cdot m^4 \cdot |S_2|^{1-1/m} &= 8 m^{1/m} \cdot m^4 \cdot (4^m m^{3m+1})^{1-1/m}\\
        &= 4 \cdot m^{1/m} \cdot m^4 \cdot (4^m m^{3m+1})^{1-1/m} + 4 \cdot 4^{m-1} m^{3m+2}\\
        &> (m-1)^{1/m} (4m^4-m+1) (4^m m^{3m+1})^{1-1/m} + (m-1) 4^m m^{3m+1}\\
        &= (m-1)^{1/m} (|M_1|-m+1) |S_2|^{1-1/m} + (m-1) |S_2|.
    \end{align*}
    Thus, by Theorem \ref{KST}, $G$ contains a copy $R_1\times R_2$ of $K_{m,m}$. Note that this $K_{m,m}$ is left-colored. Indeed, each vertex in $R_1$ colors $R_2$ with its popular color, but each vertex in $R_2$ colors $R_1$ rainbow, so the popular colors must all be distinct.\\
    
    \textbf{Case 2.}\footnote{The following argument takes an approach inspired by 
    a proof of Babai \cite{Babai}.} 
    Suppose, contrary to \textbf{Case 1}, that at least $|S_1| - 4m^4 = m$ vertices in $S_1$ use each color at most $2 m^{1/m} \cdot |S_2|^{1-1/m}$ times, and call this set $M_1'$.
    Choose $X \subseteq S_2$ which is maximal with respect to the property of $M_1' \times X$ being rainbow-colored | note that $X \neq \emptyset$ since any singleton satisfies this property by our hypothesis on $S_2$.
    Now for each $b \in S_2 \setminus X$, since $\Delta_b$ colors $M_1'$ with a rainbow pattern, by the maximality of $X$ it must be the case that there exists $a \in M_1'$ and $e \in M_1' \times X$ such that $\Delta(a,b) = \Delta(e)$.
    Moreover, the pair $(a, \Delta(e))$ can satisfy this for at most $2 m^{1/m} \cdot |S_2|^{1-1/m}$ vertices in $S_2 \setminus X$ by the definition of $M_1'$, and so
    \begin{equation*}
        |M_1' \times (M_1' \times X)| \cdot 2 m^{1/m} \cdot |S_2|^{1-1/m} \geq |S_2 \setminus X|
    \end{equation*}
    i.e.,
    \begin{equation*}
        m^2 \cdot |X| \cdot 2 m^{1/m} \cdot |S_2|^{1-1/m} \geq |S_2| - |X|.
    \end{equation*}
    Rearranging this gives 
    \begin{align*}
        |X| &\geq \frac{|S_2|}{2 m^{2+1/m} |S_2|^{1 - 1/m} + 1}\\
        &\geq \frac{|S_2|}{4 m^{2+1/m} |S_2|^{1 - 1/m}}\\
        &= \frac{|S_2|^{1/m}}{4 m^{2+1/m}},
    \end{align*}
    and the choice $|S_2| = 4^m m^{3m+1}$ guarantees that this has size at least $m$.
    Hence, the graph on $M_1' \times X$ contains a rainbow $K_{m,m}$.
\end{proof}

\section{Concluding remarks}\label{conclusion}

The most important open problem related to our theorem is, of course, to determine the order of magnitude of $\log ER(m)$.
However, if the lower bound is closer to the truth, removing the logarithmic factor in the upper bound could potentially
be a very hard problem, roughly of the same level of difficulty as determining whether the multicolor Ramsey number $R(3; q)$ of triangles
is exponential or super-exponential in the number of colors $q$.
Indeed, the best-known bounds
for $R(3; q)$
are of the form
\begin{equation*}
    2^{cq} \leq R(3;q) \leq 2^{c' q \log q},
\end{equation*}
for constants $c, c' > 0$,
and the logarithm in the exponent of this upper bound comes from the same neighborhood-chasing type argument which is responsible for it in the upper bound in (\ref{canbounds}).
As such, an improvement in either direction to Lefmann and R\"odl's bounds would undoubtedly be a major advance.

Another corollary to Rado's theorem is the following $d$-partite version.

\begin{theorem}[Rado \cite{Rado}]
    Let $m$ and $d$ be positive integers.
    Then there exists a least positive integer $n = ER^d(m)$ such that for any edge-coloring $\Delta: E(K_{n,n,\dots,n}^{(d)}) \to \N$, there exists a copy $\mathcal{K}$ of $K_{m,m,\dots,m}^{(d)}$ in $K_{n,n,\dots,n}^{(d)}$ and a subset $J \subseteq [d]$ such that for all edges $\mathbf{a} = (a_1, \dots, a_d), \mathbf{b} = (b_1, \dots, b_d) \in E(\mathcal{K})$,
    \begin{equation*}
        \Delta(\mathbf{a}) = \Delta(\mathbf{b}) \Longleftrightarrow a_i = b_i \text{ for all } i \in J.
    \end{equation*}
\end{theorem}

The current best bounds for this theorem are very far apart, with the exponential-type lower bound again coming from a random coloring and the tower-type upper bound (the height of which grows with $d$) coming from Rado's original proof.

\begin{problem}
    Give for $ER^d(m)$ either a double-exponential lower bound or a tower-type upper bound of
    height $O(1)$ as $d \to \infty$.
\end{problem}

\section{Acknowledgements}

The authors would like to thank Ayush Basu, Marcelo Sales, Cosmin Pohoata, Vojt\v{e}ch R\"odl, and Liana Yepremyan for their helpful discussions and encouragement. 


\bibliography{refs}

\begin{bibdiv}
\begin{biblist}

\bib{Abbott}{article}{
      author={Abbott, H.~L.},
       title={A note on {R}amsey's theorem},
        date={1972},
        ISSN={0008-4395,1496-4287},
     journal={Canad. Math. Bull.},
      volume={15},
       pages={9\ndash 10},
         url={https://doi.org/10.4153/CMB-1972-002-5},
      review={\MR{314673}},
}

\bib{Babai}{article}{
      author={Babai, L\'aszl\'o},
       title={An anti-{R}amsey theorem},
        date={1985},
        ISSN={0911-0119,1435-5914},
     journal={Graphs Combin.},
      volume={1},
      number={1},
       pages={23\ndash 28},
         url={https://doi.org/10.1007/BF02582925},
      review={\MR{796179}},
}

\bib{ConlonFerber}{article}{
      author={Conlon, D.},
      author={Ferber, A.},
       title={Lower bounds for multicolor {R}amsey numbers},
        date={2021},
        ISSN={0001-8708,1090-2082},
     journal={Adv. Math.},
      volume={378},
       pages={Paper No. 107528, 5},
         url={https://doi.org/10.1016/j.aim.2020.107528},
      review={\MR{4186575}},
}

\bib{ErdosRado2}{article}{
      author={Erdős, P.},
      author={Rado, R.},
       title={A combinatorial theorem},
        date={1950},
        ISSN={0024-6107,1469-7750},
     journal={J. London Math. Soc.},
      volume={25},
       pages={249\ndash 255},
         url={https://doi.org/10.1112/jlms/s1-25.4.249},
      review={\MR{37886}},
}

\bib{GMSY24}{article}{
      author={Gishboliner, L.},
      author={Milojevi{\'{c}}, A.},
      author={Sudakov, B.},
      author={Wigderson, Y.},
       title={Canonical {R}amsey numbers of sparse graphs},
        date={2024},
      eprint={2410.08644},
}

\bib{KamcevSchacht}{article}{
      author={Kam{\v{c}}ev, N.},
      author={Schacht, M.},
       title={Canonical colorings in random graphs},
        date={2023},
      eprint={2303.11206},
}

\bib{KST}{article}{
      author={K{\H{o}}vári, T.},
      author={S\'os, V.~T.},
      author={Tur\'an, P.},
       title={On a problem of {K}. {Z}arankiewicz},
        date={1954},
        ISSN={0010-1354,1730-6302},
     journal={Colloq. Math.},
      volume={3},
       pages={50\ndash 57},
         url={https://doi.org/10.4064/cm-3-1-50-57},
      review={\MR{65617}},
}

\bib{LefmannRodl}{article}{
      author={Lefmann, Hanno},
      author={R\"odl, Vojt{\v{e}}ch},
       title={On {E}rd{\H{o}}s-{R}ado numbers},
        date={1995},
        ISSN={0209-9683},
     journal={Combinatorica},
      volume={15},
      number={1},
       pages={85\ndash 104},
         url={https://doi.org/10.1007/BF01294461},
      review={\MR{1325273}},
}

\bib{Promel}{book}{
      author={Pr\"{o}mel, Hans~J\"urgen},
       title={Ramsey {T}heory for {D}iscrete {S}tructures},
   publisher={Springer, Cham},
        date={2013},
        ISBN={978-3-319-01314-5; 978-3-319-01315-2},
         url={https://doi.org/10.1007/978-3-319-01315-2},
      review={\MR{3157030}},
}

\bib{Rado}{article}{
      author={Rado, R.},
       title={Direct decomposition of partitions},
        date={1954},
        ISSN={0024-6107,1469-7750},
     journal={J. London Math. Soc.},
      volume={29},
       pages={71\ndash 83},
         url={https://doi.org/10.1112/jlms/s1-29.1.71},
      review={\MR{65616}},
}

\bib{Ramsey}{article}{
      author={Ramsey, F.~P.},
       title={On a problem of formal logic},
        date={1929},
        ISSN={0024-6115},
     journal={Proc. London Math. Soc. (2)},
      volume={30},
      number={4},
       pages={264\ndash 286},
         url={https://doi.org/10.1112/plms/s2-30.1.264},
      review={\MR{1576401}},
}

\bib{Quantitative}{article}{
      author={Reiher, Christian},
      author={R\"odl, Vojt\v~ech},
      author={Sales, Marcelo},
      author={Sames, Kevin},
      author={Schacht, Mathias},
       title={On quantitative aspects of a canonisation theorem for edge-orderings},
        date={2022},
        ISSN={0024-6107,1469-7750},
     journal={J. Lond. Math. Soc. (2)},
      volume={106},
      number={3},
       pages={2773\ndash 2803},
         url={https://doi.org/10.1112/jlms.12648},
      review={\MR{4498567}},
}

\bib{Shelah2}{article}{
      author={Shelah, Saharon},
       title={Finite canonization},
        date={1996},
        ISSN={0010-2628,1213-7243},
     journal={Comment. Math. Univ. Carolin.},
      volume={37},
      number={3},
       pages={445\ndash 456},
      review={\MR{1426909}},
}

\bib{Voigt}{article}{
      author={Voigt, Bernd},
       title={Canonizing partition theorems: diversification, products, and iterated versions},
        date={1985},
        ISSN={0097-3165,1096-0899},
     journal={J. Combin. Theory Ser. A},
      volume={40},
      number={2},
       pages={349\ndash 376},
         url={https://doi.org/10.1016/0097-3165(85)90096-2},
      review={\MR{814420}},
}

\end{biblist}
\end{bibdiv}

\end{document}